\documentclass[letterpaper,10 pt, conference]{ieeeconf}   % Comment this line out
\IEEEoverridecommandlockouts                              % This command is only
\usepackage{xcolor}
\usepackage{graphics} % for pdf, bitmapped graphics files
\usepackage{epsfig} % for postscript graphics files
\usepackage{amsmath} % assumes amsmath package installed
\usepackage{amssymb}  % assumes amsmath package installed
\usepackage[latin1]{inputenc}

\newcommand{\beq}{\begin{equation}}
\newcommand{\eeq}{\end{equation}}
\newcommand{\bea}{\begin{eqnarray}}
\newcommand{\eea}{\end{eqnarray}}
\newcommand{\bef}{\begin{figure}}
\newcommand{\eef}{\end{figure}}
\newcommand{\bsc}{\begin{scriptstyle}}
\newcommand{\esc}{\end{scriptstyle}}
\newcommand{\bd}{\begin{displaymath}}
\newcommand{\ed}{\end{displaymath}}
\newcommand{\nn}{\nonumber}

\newtheorem{lemma}{Lemma}[section]
\newtheorem{proposition}{Proposition}[section]
\newtheorem{corollary}{Corollary}[section]

\newcommand{\fig}[1]{Fig.\ \ref{#1}}
\newcommand{\alg}[1]{\begin{align} #1 \end{align}}

\def\bmat{\left[ \begin{array}}
\def\emat{\end{array} \right]}

\renewcommand{\thesubsection}{\Alph{subsection}}
\def\qed{\hfill \vrule height 7pt width 7pt depth 0pt \smallskip}

\title{ \bf Robust Kalman Filtering: Asymptotic Analysis\\ of the Least Favorable Model}

\author{Mattia~Zorzi, Bernard C.~Levy\thanks{M. Zorzi is is with the
Dipartimento di Ingegneria dell'Informazione, Universit\`a di
Padova, via Gradenigo 6/B, 35131 Padova, Italy,
({\tt\small zorzimat@dei.unipd.it})} \thanks{B. Levy is with the Department of Electrical 
and Computer Engineering, 1 Shields Avenue, University of California, Davis,
CA 95616 ({\tt\small bclevy@ucdavis.edu})}}

\begin{document}

\maketitle
\thispagestyle{empty}
\pagestyle{empty}

\begin{abstract}
We consider a robust filtering problem where the robust filter is designed according 
to the least favorable model belonging to a ball about the nominal model. In this approach,
the ball radius specifies the modeling error tolerance and the least favorable model is 
computed by performing a Riccati-like backward recursion. We show that this recursion 
converges provided that the tolerance is sufficiently small.
\end{abstract}

\section{Introduction}

Consider the problem of estimating a state process whose state-space model is known
only imperfectly. In such a situation the standard Kalman filter may perform poorly. 
Robust filtering seeks to find a state estimate which takes the model uncertainty into
account.

In this paper, we consider the robust filtering approach proposed in \cite{LN2}, see 
also \cite{LN1,HS}. The actual state-space model is assumed to belong to a ball centered 
about the nominal state-space model. The ball is formed by placing a bound on the 
Kullback-Leibler divergence between the actual and the nominal state-space model, and
the ball radius represents the modeling error tolerance. Then, the robust filter is 
designed according to the least favorable model in the ball. The resulting filter obeys 
a Kalman-like recursion which makes it very appealing for online applications \cite{ZZ}. 
Interestingly, if the ball is selected by using the $\tau$-divergence instead of the
 Kullback-Leibler divergence, the resulting filter still obeys a Kalman-like recursion 
\cite{Zor1,Zor3}. In \cite{ZL,Zor2} it was shown that when the tolerance is sufficienly 
small, the robust filter converges. Finally, it worth noting that this robust filter 
represents a generalization of risk-sensitive filters \cite{Whi,BS,LZ,HSK1} where 
large errors are severely penalized by selecting a risk-sensitivity parameter. 

It is also important to characterize the least favorable model corresponding to the 
robust filter because it can be used to evaluate the performance of an arbitrary 
filter under this least favorable situation. In \cite{LN2} it was shown that the 
least favorable model can be computed over a finite interval by first evaluating
the robust filter over the interval and then performing a backward recursion to generate the least favorable model dynamics. In this paper, we show 
that this backward recursion is a Riccati-like equation of the form 
\alg{X_{t}=A(X_{t+1}+R)^{-1}A^T+Q\nn} 
which converges provided the tolerance is sufficiently small. As a consequence, 
the least favorable model is a state-space model with constant parameters in steady 
state. The convergence of discrete-time Riccati equations with $R$ positive definite
or semi-definite has been studied in detail \cite{CM1,CM2,BGP,FN}. But in the equation 
considered here, $R$ is negative definite, and in this case only a few results are available, see \cite{Zho}. 

The outline of the paper is as follows. In Section \ref{sec:robfil} we introduce 
the robust filtering problem, in particular the backward least favorable model
recursion. In section \ref{sec:conv} we prove that the recursion converges when the 
tolerance is sufficiently small. In Section \ref{sec:perf} we show that the estimation 
error of an arbitrary stable estimator under the least favorable model is bounded. In 
Section \ref{sec:simu} some simulation results are presented. Finally, conclusions are 
presented in Section \ref{sec:conc}.

In this paper we will use the following notation. $(a,b]$ denotes an interval which is 
left-open and right-closed. Given a matrix $A\in\mathbb{R}^{n\times n}$, its spectrum 
is denoted by $\lambda(A)$ and its spectral radius is denoted by $\sigma(A)$.
We say that $A$ is (Schur) stable if $\sigma(A)<1$. $\mathcal{Q}_n$ denotes the vector 
space of symmetric matrices of dimension $n$. Given $X\in \mathcal{Q}_n$, $X>0$ 
($X\geq 0$) indicates that $X$ is positive definite (semi-definite). Given two
functions $f$ and $g$, $f(x)=o(g(x))$ around $x=\alpha$ means that $\lim_{x\rightarrow \alpha} f(x)/g(x)=0$.

\section{Robust Kalman filtering}
\label{sec:robfil}

Consider a nominal state-space model of the form
\alg{
x_{t+1} &=A x_t + B v_t \nn \\
y_t &= C x_t + D v_t  \label{nomi_model}
}
where $x_t \in \mathbb{R}^n$ is the state process, $y_t \in \mathbb{R}^p$ 
the observation process, $v_t\in\mathbb{R}^m$ is a white Gaussian noise (WGN) 
with unit variance, i.e. $\mathbb{E}[v_t v_s^T]=I \delta_{t-s}$ and $\delta_t$
denotes the Kronecker delta function. We assume that $v_t$ is independent of 
the initial state vector $x_0\sim {\cal N} (0,P_0)$, and that the pairs 
$(A,B)$ and $(A,C)$ are reachable and observable, respectively. Without loss 
of generality, we assume that $BD^T=0$. Indeed, if this is not the case we
can always rewrite (\ref{nomi_model}) with $\check A= A-ABD^T(DD^T)^{-1}C$, 
$\check B$ such that $\check B {\check B}^T=B(I-D^T(DD^T)^{-1}D)B^T$, 
$\check C=C$ and $\check D=D$. The nominal model (\ref{nomi_model}) is 
completely characterized by the transition probability density 
$\phi_t(x_{t+1},y_t|x_t)$ and by the probability density $f(x_0)$ of $x_0$. 
Let $\tilde \phi_t(x_{t+1},y_t|x_t)$ denote the transition probability density
of the actual model. We assume that the actual and nominal densities of initial 
state $x_0$ coincide, whereas $\tilde \phi_t$ belongs to a ball centered about
$\phi_t$ with radius $c>0$, hereafter called tolerance, which is specified by
\alg{
\mathcal B_t=\{ \, \tilde \phi_t\, \hbox{ s.t. } \, D_{KL}( \phi_t,\tilde \phi_t)\leq c\}.
}
Here $D_{KL}$ denotes the {\em Kullback-Leibler} divergence \cite{Kul} between 
$\phi_t$ and $\tilde \phi_t$. Note that $D_{KL} (\phi_t, \tilde \phi_t)$ is 
finite only if matrix $[\, B^T\, D^T\, ]^T$ has full row rank. Accordingly, without
loss of generality we assume that $[\, B^T\, D^T\, ]^T$ is square and invertible, 
so that $m=n+p$. Indeed it is always possible to compress the column space of this matrix 
and remove the noises which do not  affect model (\ref{nomi_model}). Let 
$Y_t=\{\, y_s,\, s\leq t\}$ and $g_t(Y_t)$ be an estimator of $x_{t+1}$ given $Y_t$. 
Adopting the minimax approach described in \cite{LN2}, a robust estimator of $x_{t+1}$ 
is obtained by solving: 
\alg{\label{minimax} 
\hat x_{t+1}=\underset{g_t\in \mathcal G_t}{\mathrm{argmin}} \max_{\tilde \phi_t \in \mathcal B_t} 
\tilde{\mathbb{E}} [\| x_{t+1}-g_t(Y_t)\|^2| Y_{t-1}] 
}
where $\tilde{\mathbb{E}}$ denotes the expectation operator taken with respect to 
the joint probability density of the actual model and $\mathcal G_t$ denotes the 
class of estimators with finite second-order moments with respect to 
$\tilde \phi_t\in \mathcal{B}_t$. In \cite{LN2}, it was shown that the robust estimator 
satisfies a Kalman-like recursion of the form
\alg{ G_t &= A V_t C^T (C V_t C^T +DD^T )^{-1}\nn\\
\hat{x}_{t+1} & = A \hat{x}_t + G_t (y_t-C \hat x_t)\nn\\
P_{t+1}&=A (V_t^{-1}+C^T(DD^T)^{-1}C)^{-1} A^T+BB^T \nn\\
V_{t+1}&=(P_{t+1}^{-1}-\theta_t I)^{-1} \label{riccati_eq}}
where $\theta_t>0$ is the unique solution to the equation $c=\gamma (P_{t+1},\theta_t)$. 
The function $\gamma$ is given by
\alg{  
\gamma(\theta,P)= \frac{1}{2}\left [\log\det(I-\theta P)+\mathrm{tr}[(I-\theta P)^{-1}]-n\right].
}
The initial conditions of the recursion are $\hat x_0=0$ and $V_0=P_0$. The least favorable 
prediction error $e_t=x_t-\hat x_t$ of the robust estimator has zero mean and covariance 
matrix $V_t$.
 
The following result is proved in \cite[Proposition 3.5]{ZL}, see also \cite{Zor2}.

\begin{proposition} \label{prop2.1}
There exists $c_{MAX}>0$ such that if $c\in (0,c_{MAX}]$, 
then for any $P_0>0$ the sequence $P_t$, $t\geq 0$, generated by (\ref{riccati_eq}) 
converges to a unique solution $P>0$, $\theta_t\rightarrow \theta$ with $\theta>0$, 
$V_t\rightarrow V$ with $V>0$ and the limit $G$ of the filtering gain $G_t$ as 
$t \rightarrow \infty$ is such that $A-GC$ is stable. Moreover, $P$ is the
unique solution of the algebraic Riccati-like equation
\alg{\label{eq:Riccati}P=A(P^{-1}-\theta I+C^T(DD^T)^{-1}C)^{-1}A^T+BB^T.}  
\end{proposition}
\smallskip\smallskip\smallskip\smallskip\smallskip\smallskip

It is possible to show that the least favorable model obtained by solving 
(\ref{minimax}) is given by \cite{LN2}
\alg{
\label{ls_model}\xi_{t+1} &= \tilde A_t \xi_t + \tilde B_t \varepsilon_t \nn \\
y_t &= \tilde C_t \xi_t + \tilde D_t \varepsilon_t  
} 
where
\alg{\label{lf_dyna} \tilde A_t &=\left[\begin{array}{cc} A & B H_t \\
0 & A-G_t C+(B-G_t D ) H_t   \end{array}\right]\nn\\
\tilde B_t & =\left[\begin{array}{c} B   \\ B-G_t D     \end{array}\right]L_t\nn \\
\tilde C_t &= \left[\begin{array}{cc}  C & D H_t    \end{array}\right],\,\, \tilde D_t=D L_t\nn\\ H_t &= \tilde K_t (B-G_t D)^T (\Omega_{t+1}^{-1}+\theta_t I)(A-G_t C)\nn\\
\tilde K_t &= [I-(B-G_t D)^T(\Omega_{t+1}^{-1}+\theta_t I ) (B-G_t D)]^{-1}
}
and $L_t$ is such that $\tilde K_t=L_t L_t^T$. In this model $\varepsilon_t$ is a
WGN with unit variance, and $\Omega_{t+1}^{-1}$ is computed by the backward recursion 
\alg{\label{back_rec} 
\Omega_{t}^{-1} = &(A-G_t C)^T[(\Omega_{t+1}^{-1}+\theta_t I)^{-1} -(B-G_t D) \times \nn\\ 
& \times (B-G_t D)^T]^{-1}(A-G_t C)
}
where if $T$ denotes the simulation horizon, the initial condition is $\Omega _{T}^{-1}=0$.

In summary, the least favorable model (\ref{ls_model}) is obtained in two steps: 

\begin{enumerate}

\item The Riccati equation (\ref{riccati_eq}) for $P_t$ is propagated forward in time
over $[0,T]$ and used to compute $G_t$ and $\theta_t$.

\item The model $(\tilde A_t,\tilde B_t,\tilde C_t,\tilde D_t)$ is obtained by propagating
(\ref{back_rec}) backward in time to evaluate $\Omega_t^{-1}$ over interval $[0,T]$.

\end{enumerate} 

It is clear that the least favorable model depends on the length $T$ of the simulation interval.
Let $\alpha,\beta$ such that $0<\alpha< \beta <1$. Then, the interval $[\alpha T,\beta T]$ is 
contained in $[0,T]$. In the next section we show that when $c>0$ is sufficient small, then 
$\Omega_t^{-1}$ converges over the interval $[\alpha T,\beta T]$ as $T$ tends to infinity.
As a consequence, the least favorable model (\ref{ls_model}) is constant over this interval.

Before establishing the convergence of the backward recursion (\ref{back_rec}), it is worth 
considering the limit case $c=0$ when the nominal and the actual models coincide. In this case,
the robust filter (\ref{riccati_eq}) reduces to the usual Kalman filter  and $\theta_t=0$ for
all $t$. Hence the limit of $\theta_t$ is $\theta=0$. By using the matrix inversion lemma, the backward recursion 
(\ref{back_rec}) with $\theta_t=0$ can be rewritten as
\alg{\Omega_{t}^{-1} = &(A-G_t C)^T[\Omega_{t+1}^{-1}-
\Omega_{t+1}^{-1}(B-G_t D)  \times \nn\\ & \times S_t  (B-G_t D)^T \Omega_{t+1}^{-1} ](A-G_t C)\nn
}
where
\alg{
S_t=[(B-G_t D)^T \Omega_{t+1}^{-1} (B-G_t D)-I]^{-1}.\nn 
}
Therefore, if $\Omega_{t+1}^{-1}=0$ then $\Omega_t^{-1}=0$. Since $\Omega_T^{-1}=0$, we conclude 
that $\Omega_t^{-1}=0$ for all $t\in[0,T]$. Accordingly, $H_t=0$ and $L_t=I$. Substituting these
expressions inside (\ref{lf_dyna}), it is then easy to verify that the least favorable model
coincides with the nominal model.

\section{Convergence of the backward recursion}
\label{sec:conv}

Suppose that the condition of Proposition \ref{prop2.1} is satisfied. Then as $t \rightarrow \infty$ the 
backward recursion (\ref{back_rec}) becomes
\alg{\label{back_rec2}
\Omega_t^{-1}=\bar A^T [(\Omega_{t+1}^{-1}+\theta I)^{-1} -\bar B \bar B^T]^{-1}\bar A
}
where the matrix $\bar A:= A-GC$ is stable, and $\bar B:= B-GD$. To ease the exposition, 
we assume that $T$ is finite and we study the convergence of (\ref{back_rec2}) as $t$ tends to
$-\infty$. This is equivalent to studying the convergence in $[\alpha T, \beta T]$ as $T$ tends to
$\infty$. Adding $\theta I $ on both sides and defining $X_t:=\Omega_t^{-1}+\theta I$ yields
the equivalent recursion 
\alg{\label{back2_rec} 
X_t=\bar  A^T(X_{t+1}^{-1}- \bar B \bar B^T)^{-1}\bar A+\theta I
}
with terminal value $X_T=\theta I$. It has the form of a Riccati equation, but an important
difference, compared to the standard case, is that in the inverse we add to $X_{t+1}^{-1}$ the
negative definite matrix $-\bar B \bar B^T$. This difference makes the convergence analysis 
nontrivial. At this point, it is useful to introduce the following map defined for $0<X<(\bar B \bar B^T)^{-1}$  
\alg{
\Theta(X):= \bar A^T (X^{-1}-\bar B \bar B^T)^{-1}\bar A+\theta I.
} 
Note that $\bar B \bar B^T$ is an invertible matrix since
\alg{
\bar B \bar B^T&= (B-GD)(B-GD)^T\nn\\ &= BB^T+ GDD^T G^T\geq BB^T
}
where $BB^T$ is invertible because $B\in \mathbb{R}^{n\times n+p}$ has full row-rank. Accordingly, 
the recursion (\ref{back2_rec}) can be rewritten as 
\alg{X_t=\Theta(X_{t+1}).}

\begin{proposition} \label{prop3.1} 
For any $0< X<(\bar B \bar B^T)^{-1}$, we have $\Theta(X)\geq \theta I$.
\end{proposition}

\proof We have 
\alg{\Theta(X)-\theta I=\bar A^T (X^{-1}-\bar B \bar B^T)^{-1}\bar A } 
where the right hand side is positive semi-definite. 
\qed 

\begin{proposition} \label{prop3.2}
The map $\Theta$ preserves the partial order of positive semi-definite matrices, so if $X_1,X_2$ 
are such that $0<X_1\leq X_2<(\bar B \bar B^T)^{-1}$, we have 
\alg{\Theta(X_1)\leq \Theta(X_2).\nn}
\end{proposition}

\proof The first variation of $\Theta(X)$ along the direction $\delta X\in \mathcal Q_n$ can be
expressed as
\alg{
\delta \Theta(X; \delta X)=& \bar A^T (X^{-1}-\bar B \bar B^T)^{-1}X^{-1}
\delta X \times \nn\\ & \times X^{-1}(X^{-1}-\bar B \bar B^T)^{-1} \bar A.
} 
Thus $\delta \Theta(X; \delta X)\geq 0$ for any $\delta X\geq 0$, so the map is nondecreasing. 
\qed

Before stating the next property of $\Theta$, we prove the following lemmas.

\begin{lemma} \label{lemma1}
It is always possible to select $c\in(0,c_{MAX}]$ such that $\theta$ is arbitrarily
small. 
\end{lemma}

\proof In \cite{ZL,Zor2} it was shown that 
\alg{ 
& \label{one}\gamma(P,\theta_1)> \gamma(P,\theta_2),\;\; \; \forall \; \theta_1>\theta_2 \hbox{ s.t. }  P\geq 0, \; P\neq 0\\\
&\gamma(P_1,\theta)\geq \gamma(P_2,\theta),\;\; \; \forall \; P_1\geq P_2\\
& \gamma(P,0)=0,\; \; \forall \; P\geq 0\\
&\label{four} \gamma(P,[0,\sigma(P)^{-1}))=[0,\infty),\; \; \forall \; P>0
} 
where (\ref{four}) means that the image of  $[0,\sigma(P)^{-1})$ under $\gamma(P,\cdot)$ 
is $[0,\infty)$. Since $c\in(0,c_{MAX}]$, by Proposition \ref{prop2.1} we have that
$P_t\rightarrow P $, $c_t\rightarrow c$, $\theta_t\rightarrow \theta$ where $c$
and $\theta$ are related by $c=\gamma(P,\theta)$. Here $P$ solves the algebraic
form of Riccati equation (\ref{riccati_eq}), so $P \geq BB^T$. In view of
(\ref{one})-(\ref{four}) it follows that $\theta\leq \tilde \theta$  where
$\tilde \theta$ is the unique solution of equation $c=\gamma(BB^T, \tilde \theta)$. 
Furthermore, the map 
\alg{\mu  \,:\, [0, \sigma(BB^T)^{-1}) &\rightarrow [0,\infty)\nn\\ 
\tilde \theta &\mapsto \gamma(BB^T, \tilde \theta)
}
is injective and continuous. Accordingly, the inverse map $\mu^{-1}\, :\, 
[0,\infty) \rightarrow [0,  \sigma(BB^T)^{-1}) $ exists and is continuous, 
in particular $\mu^{-1}(0)=0$. This means that we can always select $c>0$ such
that $\tilde \theta$ is arbitrarily small. Since $\theta\leq \tilde \theta$, 
the statement follows. \qed

It is worth noting that $\bar A$ and $\bar B$ depend on $c$ through $\theta$. Throughout 
the paper we make the following assumption.
\smallskip

{\em Assumption 1}: The map \alg{\gamma\,:\, [\,0,\; & \check \theta\,]\rightarrow \mathbb{R}^{n\times n}\times\mathbb{R}^{n\times m}\nn\\
& \theta \mapsto (\bar A,\bar B)}
is continuous for $\check \theta$ sufficiently small.

Even though Assumption 1 may appear restrictive, it holds under mild conditions on
system ($A$, $B$, $C$, $D$). Indeed, for $c\in(0,c_{MAX}]$ the unique solution of (\ref{eq:Riccati}) 
is $P=XY^{-1}$ where $[\,X^T\,Y^T ]^T$ spans the stable deflating subspace of 
regular matrix pencil $sL -M$ \cite{PLS}, where
\alg{
L=\left[\begin{array}{cc}A^T & 0  \\ -BB^T & I \end{array}\right],\;\; M=\left[\begin{array}{cc}I &
C^T(DD^T)^{-1}C-\theta I  \\0 & A \end{array}\right]. \nn
}
Conditions for the continuity of such subspaces are given in \cite{DK}. Accordingly, the
map $\theta \mapsto P$ is continuous over $[\,0\,\,\check \theta]$ with $\check \theta$ 
small enough. Since the map $P\mapsto (\bar A,\bar B)$ is continuous, we conclude that 
$\gamma$ is continuous for $\check \theta$ sufficiently small.

\begin{lemma} \label{lemma2} For $c\in(0,c_{MAX}]$ sufficiently small, there exists 
$\rho\in (1, \sigma(\bar A)^{-1})$ such that 
\alg{
\label{cond_lyap}(1-\rho^{-2})\Sigma_q^{-1}-\bar B \bar B^T \geq 0
}
where $\Sigma_\rho$ is the unique solution of the algebraic Lyapunov equation (ALE)
\alg{
\label{eq_lyap}\Sigma_\rho =\rho^2 \bar A^T \Sigma_\rho \bar A +\theta I.
}
\end{lemma}

\proof First, note that $\rho \bar A $ is a stable matrix. Then, the solution of (\ref{eq_lyap}) 
is given by
\alg{\label{expl_sol} 
\Sigma_\rho=\theta \sum_{k\geq 0} \rho^{2k}( \bar A^T)^k \bar A^k
} 
which is positive definite. Note that
\alg{\Sigma_\rho\leq \theta \sum_{k\geq 0} \rho^{2k}\sigma(\bar A)^{2k}I=\frac{\theta}{1-\rho^2 \sigma(\bar A)^2}I\nn}
and thus $\Sigma_\rho^{-1}\geq (1-\rho^2 \sigma(\bar A)^2)/\theta I$.
In view of Assumption 1, for $\theta$ sufficiently small we have 
\alg{\sigma(\bar A)^2=\sigma(\bar A_0)^2 +o(1)\nn} 
where $\bar A_0=A-G_0C$, $G_0=AP^{(0)}C^T(CP^{(0)}C^T+DD^T)^{-1}$ and $P^{(0)}$ is 
the unique solution of (\ref{eq:Riccati}) with $\theta=0$. As a consequence,
\alg{
(1-\rho^{-2})\Sigma^{-1}_\rho\geq(1-\rho^{-2}) \frac{1-\rho^2 (\sigma(\bar A_0)+o(1))^2}{\theta} I
. \label{lowbound} 
} 
We can always choose $\rho$ in the range $(1, \sigma(\bar{A}_0)^{-1})$ such that $(1-\rho^{-2}) 
(1-\rho^2 \sigma^2 (\bar{A}_0))$ is positive. By Lemma \ref{lemma1}. we can also select
$c\in[0,c_{MAX}]$ sufficiently small so that $\theta$ is small enough that the 
scaled identity matrix on the right hand side of (\ref{lowbound}) upper bounds
$\bar{B}\bar{B}^T$.
\qed

Let $\bar{c} \in (0,c_{MAX}]$ be a value of $c$ such that Lemma \ref{lemma2} is satisfied,
so that (\ref{cond_lyap}) holds for a certain $\rho$ and $\theta$. Then it is useful to
observe that for any $c\in(0,\bar c)$, the equation (\ref{cond_lyap}) still holds with 
the same value for $\rho$ but with a smaller value for $\theta$. 

\begin{corollary}
For any $c\in(0,\bar c]$, we have $\Sigma_\rho< (\bar B \bar B^T)^{-1}$.
\end{corollary}

\proof Since (\ref{cond_lyap}) holds for a suitable $\rho>0$, we have 
\alg{
\Sigma^{-1}_\rho \geq \rho^{-2} \Sigma^{-1}_\rho+\bar B \bar B^T > \bar B \bar B^T\nn
}
which implies $\Sigma_\rho< (\bar B \bar B^T)^{-1}$.
\qed

We are now ready to state the third property of the map $\Theta$.

\begin{proposition} \label{prop3.3}
Consider the compact set 
\alg{
\mathcal  C=\{ X\in \mathcal Q_n \hbox{ s.t. } \theta I \leq X\leq \Sigma_\rho\}\nn
} 
where $\Sigma_\rho$ is computed as in Lemma \ref{lemma2}. If $c\in(0,\bar c]$ then 
$\Theta(X)\in \mathcal C$ for any $X\in \mathcal C$. 
\end{proposition}

\proof 
First, observe that $\mathcal C$  is a nonempty set. Indeed, by (\ref{expl_sol}) we have 
$\Sigma_\rho \geq \theta I$, so that $\theta I \in \mathcal C$. Since $c \in (0,\bar c]$, 
by Lemma \ref{lemma2} the inequality (\ref{cond_lyap}) holds for some $\rho\in (1, \sigma(\bar A)^{-1})$, 
and thus
\alg{
&\Sigma_\rho^{-1}- \bar B \bar B^T \geq \rho^{-2}\Sigma_\rho^{-1}\nn\\
&(\Sigma_\rho^{-1}- \bar B \bar B^T )^{-1}\leq \rho^{2}\Sigma_\rho\nn\\
&\bar A^T (\Sigma_\rho^{-1}- \bar B \bar B^T )^{-1}\bar A +\theta I\leq \rho^{2}\bar A^T \Sigma_\rho \bar A +\theta I\nn\\
& \label{upper_theta}\Theta(\Sigma_\rho)\leq \Sigma_\rho.
} 
Assume that $X\in \mathcal C$. Since $X\leq \Sigma_\rho$, the nondecreasing property
of $\Theta$ and (\ref{upper_theta}) imply
\alg{
\Theta(X)\leq \Theta (\Sigma_\rho) \leq \Sigma_\rho . \nn
}
Since $X\geq \theta I$, we have
\alg{
\Theta (X) \geq \Theta(\theta I) \geq \theta I }
where we exploited again the nondecreasing property of $\Theta$ and Proposition \ref{prop3.1}. 
We conclude that $\Theta(X)\in \mathcal C$.
\qed

\begin{proposition}
Consider the sequence $X_t$ satisfying the backward recursion 
\alg{
\label{iter_fin}X_t=\Theta(X_{t+1}), \; \; \; X_T=\theta I.
}
For $c\in(0,\bar c]$, the sequence belongs to $\mathcal C$ and is nondecreasing. 
Thus as $t \rightarrow -\infty$, $X_t$ converges to $X\in \mathcal C$ which is a
solution of the algebraic Riccati equation 
\alg{ \label{ric_alg}  
X=\bar A^T(X^{-1}-\bar B \bar B^T)^{-1} \bar A+\theta I.
}
\end{proposition}

\proof We prove the first two statements by induction. We start by showing that $X_t\in \mathcal C$ 
for any $t$. We know that $X_T\in \mathcal C$ because $\mathcal C$ contains $\theta I$. Assume 
that $X_{t+1}\in \mathcal C$, then Proposition \ref{prop3.3} implies that $X_t=\Theta(X_{t+1})
\in \mathcal C$. This proves the first claim. 

Next we show that the sequence is nondecreasing. We observe that 
\alg{
X_{T-1}=\Theta(X_T)=\Theta(\theta I) \geq \theta I =X_T
}
where we exploited the nondecreasing property of $\Theta$, see Propositions \ref{prop3.2} and 
\ref{prop3.1}. Assume that $X_t\geq X_{t+1}$, then 
\alg{
X_{t-1}=\Theta(X_t)\geq \Theta(X_{t+1})=X_t,
}
so by induction the sequence is nondecreasing.   

The convergence follows from the fact that the sequence is nondecreasing and belongs to a 
compact set. \qed 

Since $X_t=\Omega_t^{-1}+\theta I$, we have the following result.
\smallskip\smallskip

\begin{corollary} \label{cor_omega}
For $c\in(0,\bar c]$, the sequence $\Omega_t^{-1}$ generated by (\ref{back_rec2}) converges 
to $\Omega^{-1}$ as $t\rightarrow -\infty$ where $\Omega^{-1}$ is such that $0 \leq \Omega^{-1}
\leq \Sigma_\rho -\theta I$ for some $\rho \in (1, \sigma(\bar A)^{-1})$ satisfying (\ref{cond_lyap}). 
Furthermore
\alg{
& H_t\rightarrow H,\; \; \tilde K_t\rightarrow \tilde K,\; \;L_t\rightarrow L\nn\\
& \tilde A_t\rightarrow  \tilde A,\; \; \tilde B_t\rightarrow  \tilde B\nn\\
&  \tilde C_t\rightarrow  \tilde C,\; \; \tilde D_t\rightarrow  \tilde D.
}
\end{corollary}
\smallskip\smallskip\smallskip\smallskip\smallskip\smallskip
It is worth noting that the algebraic equation (\ref{ric_alg}) may admit several
positive definite solutions. Indeed, in the scalar case, equation (\ref{ric_alg}) becomes
\alg{\label{ric_alg_scal} 
x=\frac{\bar a^2}{x^{-1}-\bar b^2}+\theta
}
or equivalently
\alg{
\bar b^2 x^2-(1-\bar a^2+\bar b^2 \theta)x+\theta=0.\nn
} 
For small $\theta >0$, the discriminant of this equation is positive, so the equation 
has two positive real solutions since the coefficient $1 -\bar{a}^2 - \bar{b}^2 \theta$ is
positive. For $\bar a=0.1$, $\bar b=1$ and $\theta=0.1$ we obtain the two solutions 
$x_1\approx 0.99$ and $x_2\approx 0.10$. It is not difficult to see that (\ref{ric_alg_scal}) 
can be rewritten as a Lyapunov equation 
\alg{
\label{alg_ly}x=(\bar a-j\bar b)^2x+\bar b^2 -j^2
}
where $j=\bar a x \bar b /(\bar b^2 x-1)$. Let $f:=\bar a -j \bar b$ be the ``feedback'' matrix
and $f_1,f_2$ denote the values corresponding to $x_1$ and $x_2$, respectively. Then we have 
$f_1\approx8.9$ and $f_2\approx 0.11$. In view of (\ref{alg_ly}), this means that $x_1$ is 
a stabilizing solution of (\ref{back2_rec}) whereas $x_2$ corresponds to an unstable one.
Accordingly, the limit of the sequence (\ref{iter_fin}) is $x_2$. In the general case (i.e., for $n>1$) 
the algebraic Riccati equation (\ref{ric_alg}) can be rewritten as 
\alg{
X=(\bar A-\bar B J^T)^T X(\bar A-\bar B J^T) +\bar B \bar B^T-JJ^T\nn
}
where $J= \bar A ^T X \bar B (\bar B X \bar B^T- I)^{-1}$. However, the reasoning used in 
the scalar case cannot be applied since the matrix $\bar B \bar B^T-JJ^T$ is indefinite.

\begin{proposition} \label{prop3.5}
For $c\in(0,\bar c]$ sufficiently small, the limit $X$ of (\ref{iter_fin}) is a stabilizing 
solution of (\ref{ric_alg}) in the sense that the matrix $\bar A^T- J \bar B^T$ is stable.
\end{proposition}

\proof Let $X_\theta$ be the limit of the sequence in (\ref{iter_fin}) where we made explicit 
its dependence on $\theta$. Notice that $\rho$ does not depend on $\theta$. Indeed, if a 
certain $\rho$ satisfies (\ref{cond_lyap}) for a given $\theta$, then the same $\rho$ 
satisfies (\ref{cond_lyap}) with $\theta^\prime$ such that $0<\theta^\prime \leq \theta$. 
Since $X_\theta\in \mathcal C$, we have that $\theta I \leq X_\theta \leq 
\theta \sum_{k\geq 0} \rho^{2k} (\bar A^T)^k\bar A^k$. Let $Q_{\theta}$ be such that 
$X_\theta=\theta Q_\theta$. Hence $Q_\theta\geq I$. Observe that
\alg{ 
M_\theta& := \bar A^T- J \bar B^T\nn\\ 
&= \bar A^T [X_\theta-X_\theta \bar B (\bar B X_\theta \bar B^T -I)^{-1} \bar B^T X_\theta] X_\theta^{-1}\nn\\
&= \bar A^T (X_\theta^{-1}-\bar B \bar B^T)^{-1}X_\theta^{-1}\nn\\
&= \bar A^T (\theta^{-1}Q^{-1}_\theta-\bar B \bar B^T)^{-1}\theta^{-1}Q_\theta^{-1}\nn\\
&= \bar A^T (Q_\theta^{-1}-\theta \bar B \bar B^T)^{-1}Q^{-1}_\theta . \label{mtheta}
}
For $\theta$ sufficiently small, by Assumption 1 we have $\bar B \bar B^T=\bar B_0 \bar B^T_0+o(1)$ 
where $\bar B_0=B-G_0D$ and $G_0$ has been defined in the proof of Lemma  \ref{lemma2}. Accordingly,
\alg{(
Q_\theta^{-1}-\theta \bar B\bar B^T)^{-1}=Q_\theta+o(1)  ,
}
which after substitution inside (\ref{mtheta}) gives
\alg{M_\theta = \bar A^T + o(1) .
}
The map $\theta \mapsto \lambda(M_\theta)$ is a continuous function for $\theta> 0$ since 
the mapping from the entries of a matrix to its spectrum is continuous. Hence for $\theta $ 
sufficiently small, the matrix $M_\theta$ is stable. By Lemma \ref{lemma1} we conclude that 
if we select $c\in (0,c_{MAX}]$ sufficiently small, the matrix $M_\theta$ will be stable. 
\qed 

\section{Performance analysis} 
\label{sec:perf}

We want to evaluate the performance of an arbitrary estimator 
\alg{ \label{arb_filter}
\hat x_{t+1}^\prime=A \hat x_t^\prime +G_t^{\prime} (y_t-C \hat x_t^\prime)
 }
under the least favorable model (\ref{ls_model}) in steady state, i.e. with $\tilde A_t$, 
$\tilde B_t$, $\tilde C_t$  and $\tilde D_t$ constant. Note that the steady state condition 
is guaranteed under the assumption that 
$c\in(0,\bar c]$. Recall that $e_t$ denotes the least favorable prediction error of the robust 
filter (\ref{riccati_eq}). Let $e_t^\prime=x_t -\hat x_t^\prime$ be the prediction error of 
filter (\ref{arb_filter}).  Let $\mathbf{e}_t=[\,e^{\prime\,T}_t \, e_t^T\,]^T$. In \cite{LN2}
it was shown that the dynamics of $\mathbf{e}_t$ are given by
 \alg{\label{e_ls}\mathbf e_{t+1}=F_t \mathbf{e}_t+M_t \varepsilon_t}
where 
\alg{
F_t:=& \tilde A-\left[\begin{array}{c} G_t^\prime  \\  0\end{array}\right]
\tilde C ,\; \;
 M_t:=& \tilde B-\left[\begin{array}{c} G_t^\prime  \\  0\end{array}\right]
\tilde D\nn
}
and $ \varepsilon_t$ is a WGN with unit variance. Then the covariance matrix $\Pi_t$ of $\mathbf e_t$
obeys the Lyapunov equation
\alg{
\label{lyap_rec}\Pi_{t+1}=& F_t\Pi_t F_t^T+M_t M_t^T
}
with initial condition $\Pi_0=I_2\otimes V_0$.  

From (\ref{e_ls}) it is clear that the mean of the prediction error $e_t^\prime$ is zero. 
Next, we show that the covariance matrix of $e^\prime_t$ converges to a constant matrix 
and is bounded provided that $c$ is sufficiently small. To do so, we use the following 
result \cite[Theorem 1]{CS}.

\begin{lemma} \label{prop_sayed}
Consider the time-varying Lyapunov equation
\alg{
Y_{t+1}=\mathcal F_t Y_t \mathcal  F_t^T+\mathcal  R_t\nn
}
where $\mathcal F_t$ and $\mathcal R_t$ converges to $\mathcal F$ and $\mathcal  R$, respectively,
as $t\rightarrow \infty$ with $\mathcal F$ stable. Then $Y_t$ converges to the unique 
solution $Y$ of the Lyapunov equation:
\alg{
Y=\mathcal F Y \mathcal  F^T+\mathcal R.\nn
}
\end{lemma}
\smallskip\smallskip

\begin{proposition} \label{prop4.1} Assume that the gain $G_t^\prime$ in (\ref{arb_filter})
converges to a matrix $G^\prime$ such that $A-G^\prime C$ is stable. Then, for $c\in(0, \bar c]$ 
sufficiently small the recursion (\ref{lyap_rec}) converges to the solution $\Pi$ of the
Lyapunov equation
\alg{
\Pi=& F\Pi F^T+MM^T\nn
}
where 
\alg{  F:= \tilde A- \left[\begin{array}{c} G^\prime \\ 0 \end{array}\right] \tilde C,\; \;M:= \left(  \tilde B-\left[\begin{array}{c} G^\prime  \\  0\end{array}\right]
\tilde D\right).\nn
}
\end{proposition}
\vspace{3mm}
\proof First, we prove that the  matrix 
\alg{
F= \left[\begin{array}{cc}  A- G^\prime  C & (B-G^{\prime}D)H \\
0 & A-GC+(B-GD)H 
\end{array}\right]
}  
is stable. Since $F$ is an upper block-triangular matrix, it is sufficient to show that 
its two diagonal blocks are stable. The matrix $A- G^{\prime} C$ is stable by assumption. 
Next, by recalling that $\bar A= A-GC$, $\bar B=B-GD$, $H= \tilde K\bar B^T X \bar A$ and
$\tilde K=(I- \bar B^T X \bar B)^{-1}$, the (2,2) block of $F$ can be expressed as
\alg{
\bar A&+\bar B (I- \bar B X \bar B^T)^{-1}\bar B^T X \bar A\nn\\
&= \bar A- \bar B(\bar B X \bar B^T-I)^{-1}\bar B^T X \bar A\nn\\
&= \bar A- \bar B J^T} 
which has the same eigenvalues of $\bar A^T -J \bar B^T$. By Proposition \ref{prop3.5}, 
this matrix is stable provided that $c$ is sufficiently small. The conditions of Proposition
\ref{prop_sayed} are satisfied since $F_t$ converges to $F$ with $F$ stable and $M_t$ 
converges to $M$, and thus $M_tM_t^T$ converges to $MM^T$ as $t\rightarrow \infty$. 
Hence $\Pi_t$ converges to $\Pi$. 
\qed\\
  
\begin{corollary} 
Under the assumption that $c\in(0,\bar{c}]$ is sufficiently small, the prediction error $e_t^\prime$ 
of the filter (\ref{arb_filter}) under the least favorable model (in steady state) has zero mean and 
bounded variance. 
\end{corollary} 

\section{Simulation Example}
\label{sec:simu}

Consider the state-space model 
\alg{ \label{example_model}
  A&=\left[
       \begin{array}{cc}
         0.1 & 1 \\
         0 & 1.2 \\
       \end{array}
     \right], \;  B=0.01I_2\nn\\
   C&=\left[
         \begin{array}{cc}
           1 & -1 \\
         \end{array}
       \right],
       \; D=0.04 .
} 
Note that the pairs $(A,B)$ and $(A,C)$ are reachable and observable, respectively.
Using the procedure of \cite[Proposition 3.5]{ZL}, it results that the robust filter
(\ref{riccati_eq}) converges for $c\in(0,c_{MAX}]$, with $c_{MAX} =0.1879$. 

\begin{figure}
\includegraphics[width=\columnwidth]{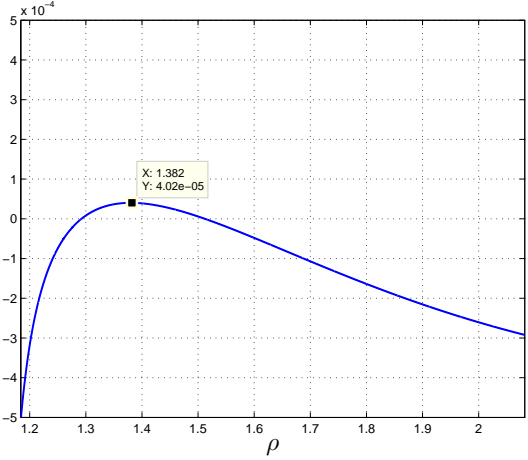}
\caption{Minimum eigenvalue of matrix $(1-\rho^{-2}) \Sigma_\rho^{-1}-\bar B \bar B^T $ as 
a function of $\rho$ for $c= c_{MAX}$.}
\label{fig_check}
\end{figure}

\begin{figure}
\includegraphics[width=\columnwidth]{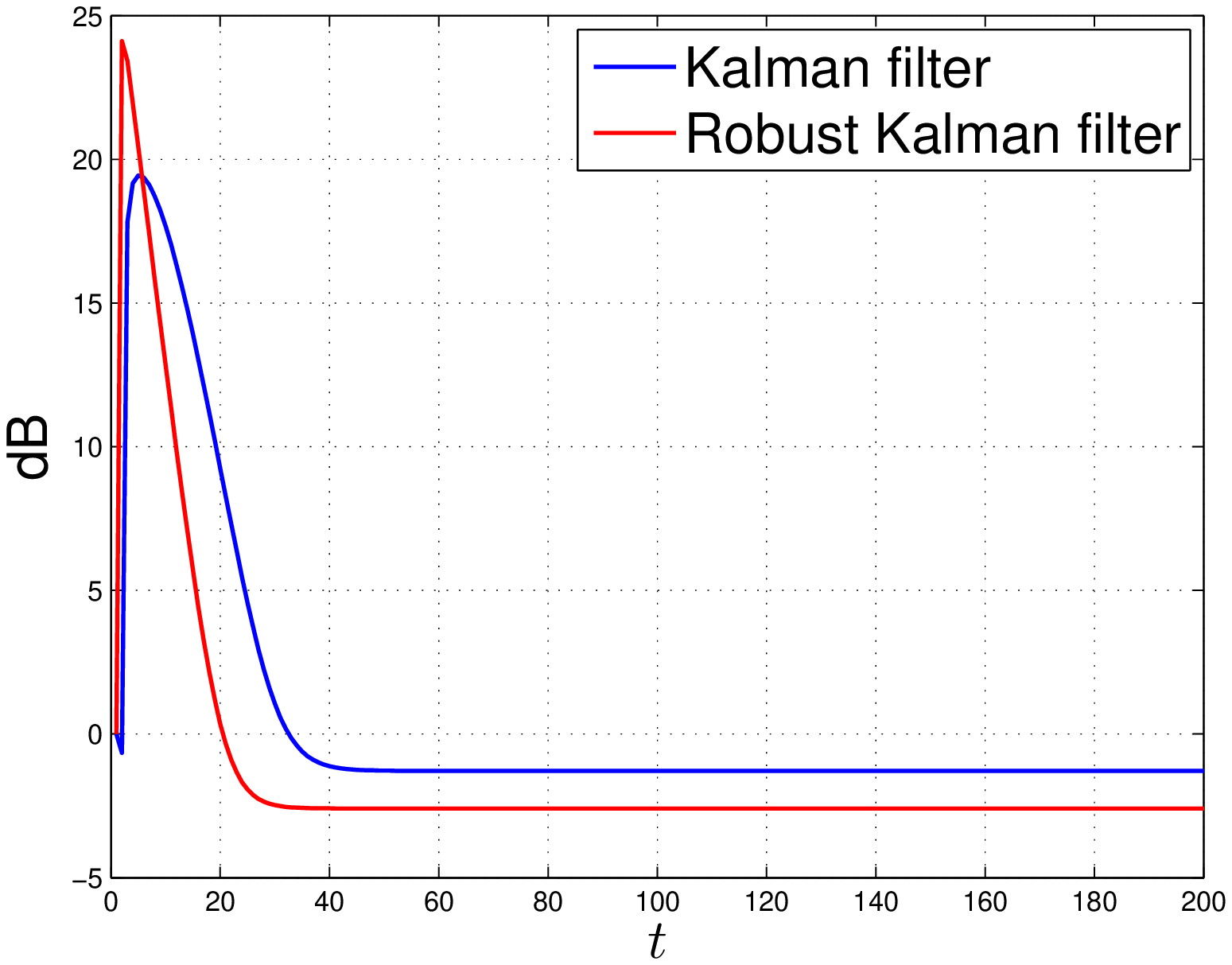}
\caption{Variance (in decibel) of the first component of the prediction error of the 
Kalman and robust filters for the least favorable model.} 
\label{fig_e1}
\end{figure}

\begin{figure}
\includegraphics[width=\columnwidth]{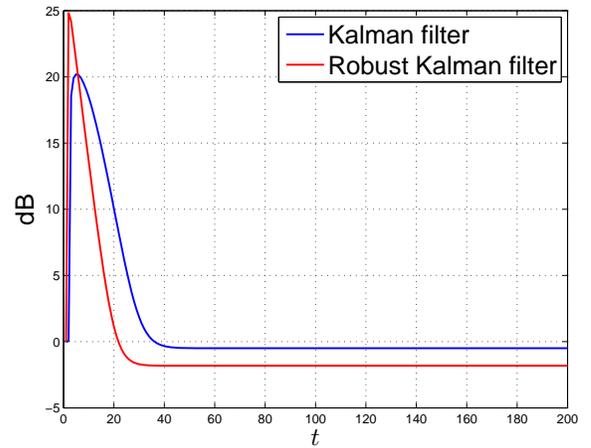}
\caption{Variance (in decibel) of the second component of the prediction error of the 
Kalman and robust filters for the least favorable model.}
\label{fig_e2}
\end{figure}

The minimum eigenvalue of $(1-\rho^{-2}) \Sigma_\rho^{-1}-\bar B \bar B^T $ is depicted 
in \fig{fig_check} as a function of $\rho$ for $c = c_{MAX}$. 
For $c=c_{MAX}$, we see that when $\rho=1.382$, the minimum eigenvalue is $4.02\cdot 10^{-5}$, 
so the matrix is positive definite and $\bar c=c_{MAX}$. Consider the 
sequence generated by (\ref{back_rec}) for $c=c_{MAX}$. We have
\alg{
\Sigma_\rho \approx10^2 \cdot\left[
       \begin{array}{cc}
         5.89&  -5.03  \\
 -5.03  & 4.31    \\
       \end{array}
     \right]  .
 \nn}
and iteration  (\ref{back_rec}) converges to 
 \alg{\Omega^{-1}\approx10^2 \cdot\left[
       \begin{array}{cc}
      4.56 & -3.90      \\
   -3.90 &  3.34  \\
       \end{array}
     \right]. \nn
} 
Furthermore, the matrix $\bar A^T- J\bar B^T$ has for eigenvalues  
$0.8373$, $0.0892$, so it is stable. Finally, Figures \ref{fig_e1} and \ref{fig_e2} 
depict the variances of the first and second component of prediction error of the Kalman 
filter and robust filter Kalman for the steady-state least favorable model. As expected,
both variances converge to a constant value and for both components, the performance of
the robust filter is approximately 1.5 dB lower than that of the Kalman filter.

\section{Conclusion} 
\label{sec:conc}

We have considered a robust filtering problem, where the minimum variance estimator is 
designed according to the least favorable model belonging to a ball about the nominal model
and with a certain radius corresponding to the modeling tolerance. We showed that as long as
the model tolerance does not exceed a maximum value $\bar{c}$, the least favorable model
converges to a constant model. Furthermore, as long as the tolerance is sufficiently small,
the covariance matrix of the prediction error for any stable filter remains bounded 
when applied to the steady-state least favorable model.

\end{document}